\documentclass[11pt]{article}
\usepackage{amssymb,amsmath,amsthm}
\usepackage[nobysame]{amsrefs}
\usepackage{verbatim}
\usepackage[top=3cm, bottom=3cm, left=2.5cm, right=2.5cm]{geometry}
\usepackage{tikz}
\usepackage{caption}
\usepackage{subcaption}

\newtheorem{theorem}{Theorem}[section]
\newtheorem{lemma}[theorem]{Lemma}
\newtheorem{proposition}[theorem]{Proposition}
\newtheorem{corollary}[theorem]{Corollary}

\title{The saturation number of carbon nanocones and nanotubes}
\author{Taylor Short\thanks{This work was partially supported by a SPARC Graduate Research Grant from the Office of the Vice President for Research at the University of South Carolina.}\\Department of Mathematics\\Grand Valley State University\\shorttay@gvsu.edu}
\date{\today}

\begin{document}
%%%%%%%%%%%%%%%%%%%%%%%%%%%%%%%%%%%%%%%%%%%%%%%%%%%%
\maketitle

\begin{abstract}
The saturation number of a graph is the cardinality of a smallest maximal matching. This paper presents bounds for the saturation number of carbon nanocones which are asymptotically equal. The same techniques are applied for the saturation number of certain families of carbon nanotubes, which improve previous results and in one case, yields the exact value.

\end{abstract}

\section{Introduction}

Throughout this paper $G$ is an $n$-vertex, simple graph with vertex set $V(G)$ and edge set $E(G)$. A \textit{matching} $M$ in a graph $G$ is a collection of edges of $G$ such that no two edges from $M$ share a vertex. The cardinality of $M$ is called the \textit{size} of the matching. A matching $M$ is a \textit{maximum matching} if there is no matching in $G$ with greater size. The \textit{matching number} $\nu (G)$ of $G$ is the cardinality of any maximum matching in $G$. Since each vertex can be incident to at most one edge of a matching, it follows that $\nu (G) \le \lfloor n/2 \rfloor$ for any graph $G$. If every vertex of $G$ is incident with an edge in $M$, then $M$ is called a \textit{perfect} matching and such graphs have $\nu (G)=n/2$. It is clear that perfect matchings are also maximum matchings but the converse is not generally true.

Matchings serve as models of many phenomena across the sciences. An important motivation for their study arose from chemistry, when it was observed that the stability of benzenoid compounds is related to the number of perfect matchings, also known as \textit{Kekul\'{e} structures}, in the corresponding chemical graphs. For a survey of these results, see \cite{gutman1}. With the discovery of fullerenes in 1985 \cite{kroto1}, the desire to identify properties characteristic for stable fullerenes led to the enumeration of perfect matchings \cite{doslic4,zhang1,kardos1,doslic3} in these corresponding graphs. 

Maximum matchings give one way to quantify the largeness of a matching. Both the enumerative and structural properties of maximum matchings are well studied and well understood, see \cite{lovasz1} for a general background on such matching theory.

There is yet another way to quantify the largeness of a matching. A \textit{maximal matching} in $G$ is a matching that cannot be extended to a larger matching in $G$. Clearly, every maximum matching is also maximal but the opposite is usually not true. Chemically, maximal matchings model the adsorption of dimers to a molecule, where each dimer bonds to a pair of adjacent atoms in the molecule. Any such adsorption pattern corresponds to a matching in the graph of the molecule, and once no further adsorption is possible, such a matching must be maximal. The best case of adsorption can be viewed as a maximum matching, while the worst case concerns the smallest possible maximal matching. This idea gives rise to the study of the \textit{saturation number} of a graph $G$, which is the cardinality of any smallest maximal matching in $G$. Thus the saturation number is a measure of how inefficient the adsorption process can be. Aside from chemistry, the saturation number has a number of interesting applications in networks, engineering, etc. The saturation number of a graph is equal to the cardinality of an independent edge dominating set. Finding an independent edge dominating set in a graph is an NP-Hard problem \cite{Yannakakis1980}.

Maximal matchings are much less understood than their maximum counterparts. Some work has been done on enumerating maximal matchings in certain chemical graphs \cite{doslic5,short1} but this area remains largely unexplored. Structural properties, such as the saturation number, have been studied for benzenoid graphs \cite{doslic2}, fullerenes \cite{andova,doslic7,doslic6}, and nanotubes \cite{tratnik1}. The paper \cite{doslic2} mentions the saturation number of nanocones as an interesting, unexplored avenue of study.

This paper considers both nanocones and nanotubes, which are carbon networks situated between graphene and fullerene in terms of structure. New upper and lower bounds on the saturation number of nanocones are established, which are asymptotically equal. In addition, lower bounds for the saturation number of two classes of nanotubes are presented, which improve recent results \cite{tratnik1}.

\section{Statement of Results}

A \textit{hexagonal patch}, or \textit{patch} for short, is a planar graph where all faces are hexagons except one \textit{outer} or \textit{boundary} face. All internal vertices have degree 3 and all vertices on the outer face have degree 2 or 3. For a face $F$ in a planar graph $G$, let $n_2(F)$ be the number of degree 2 vertices incident to $F$ and let $n_2=n_2(G)$ be the total number of degree 2 vertices in a graph $G$. 

Next it will be useful to introduce some definitions utilized in \cite{graves1,graves2,graver1,brinkmann1,bornhoft1}. The \textit{boundary code} of a patch is described by a sequence of 2's and 3's corresponding to the degree of the vertices on the boundary of the patch in cyclic order. A \textit{break edge} is an edge on the boundary whose endpoints are both degree 2. A \textit{bend edge} is an edge on the boundary whose endpoints are both degree 3. 

This paper limits itself to patches with nice boundaries. A patch is \textit{pseudoconvex} if it does not contain any bend edges. A \textit{side} of a patch is a path on the boundary between a consecutive pair of break edges, including the break edges, and the \textit{length} of a side is the number of degree 3 vertices on the side.

A \textit{defect} in a patch is a non-hexagonal face. A defect is \textit{internal} if all vertices incident to the face are degree 3. A defect is \textit{external} if there are degree 2 vertices incident to the face. Using this terminology, the outer face of a patch can also be called an external defect. Patches can have more than one external defect. In such a graph, any face incident to degree 2 vertices could the outer face in a planar drawing.

The following theorem is the main tool in proving lower bounds on the saturation number. This theorem is a generalization of the theorem proven for fullerenes in \cite{andova}, in that a fullerene graph can be viewed as a patch containing exactly 12 pentagonal defects and no vertices of degree 2 (i.e. fullerenes have no external defects). For the sake of consistency, the proof in Section \ref{seclowerbound} uses similar terminology and structure to what was presented in \cite{andova}.

\begin{theorem} \label{lowerbound}
Let $G$ be a pseudoconvex patch with $n$ vertices, $o_k$ internal defects which are $k$-gonal where $k$ is odd, $e_k$ internal defects which are $k$-gonal where $k\neq 6$ is even, and $n_2$ vertices of degree 2. Then
$$
s(G) \ge \frac{n}{3}-\frac{1}{18}\left( \sum _{k \text{ odd}}(k-2)o_k+ \sum _{k \text{ even}} ke_k \right)-\frac{n_2}{6}.
$$
\end{theorem}

\subsection{Nanocones}

Generally speaking, nanocones are planar graphs where the majority of faces are hexagons, along with some non-hexagonal faces, most commonly pentagons, in addition to the outer face. This paper considers these pentagonal defect nanocones as well as nanocones with a single $k$-gonal defect.

A \textit{single-defect $k$-gonal nanocone with $\ell$ layers}, $CNC_k(\ell)$, is obtained by taking a cycle on $k \ge 3$ vertices, $C_k$, and surrounding it with $\ell$ concentric layers of hexagons. Using previous terminology, a single-defect $k$-gonal nanocone is a pseudoconvex patch with a single $k$-gonal defect at its apex. By induction, it follows that there are $k{\ell+1\choose 2}$ hexagonal faces, $k(\ell+1)^2$ total vertices, and $(2\ell+1)k$ external vertices. There are $k\ell$ external vertices of degree 3 and $k(\ell+1)$ vertices of degree 2. The following Corollary is an immediate consequence of Theorem \ref{lowerbound}.

\begin{corollary} \label{lowersd}
$$
s(CNC_k(\ell)) \ge \begin{cases} \frac{k(\ell+1)^2}{3}-\frac{k-2}{18} - \frac{k(\ell+1)}{6} & \text{if } k \text{ is odd} \\ \frac{k(\ell+1)^2}{3}-\frac{k}{18} - \frac{k(\ell+1)}{6} & \text{if } k \text{ is even}\end{cases}
$$
\end{corollary}

\begin{figure}[h]
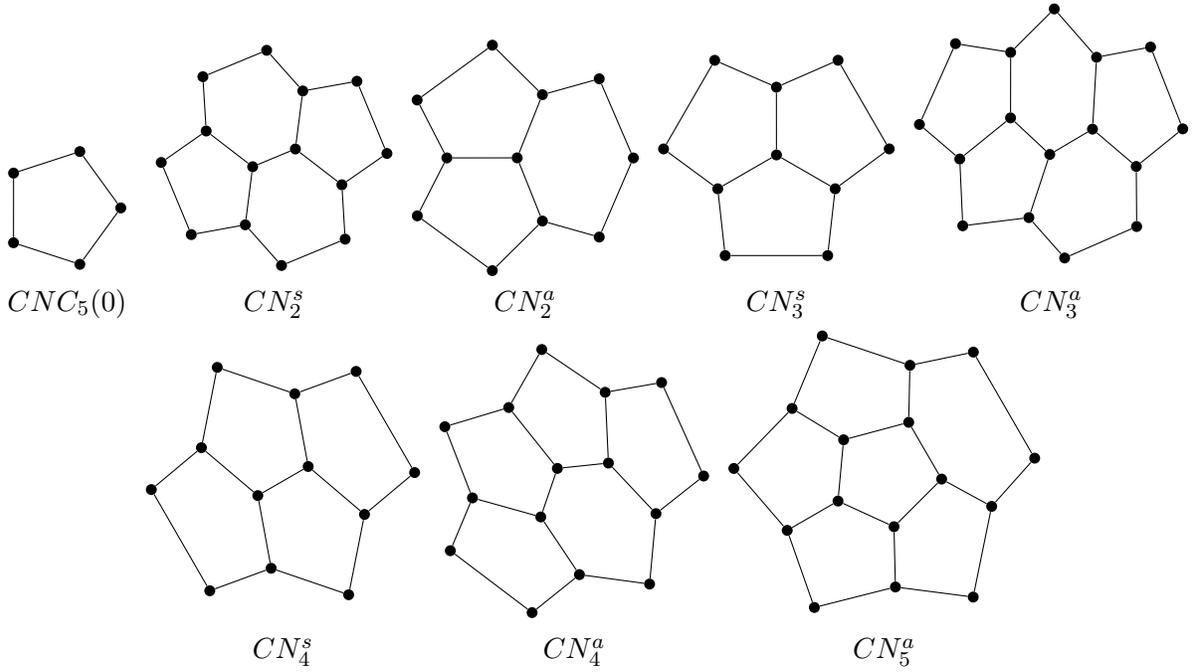

\begin{center}
  \input{cone1} \input{cone2s.tex} \input{cone2n} \input{cone3s} \input{cone3n} \input{cone4s} \input{cone4n} \input{cone5s}
\end{center}
  \caption{The 8 configurations of pentagons and or hexagons for pentagonal defect nanocones.} \label{cones}
\end{figure}

\textit{Pentagonal defect nanocones} are pseudoconvex patches with $p$ pentagons, where $1\le p \le 5$. While many arrangements of pentagons and hexagons are possible, a classification result first from \cite{klein1,klein2} and then independently in \cite{justus1} shows that it suffices to consider the 8 configurations of pentagons and or hexagons in Figure \ref{cones}. Note that the single pentagon configuration is merely $CNC_5(0)$. 

A \textit{pentagonal defect nanocone with $i$ pentagons and $\ell$ layers}, $CN_i^j(\ell)$, $i\in \{2,3,4,5\}$ and $j \in \{s,a\}$, is defined to be the configuration $CN_i^j$ in Figure \ref{cones} surrounded by $\ell$ concentric layers of hexagons. For reference, the use of $s$ in the superscript designates a symmetric configuration as drawn in Figure \ref{cones} and the use of $a$ represents asymmetric. The configurations $CNC_5(0)$ and $CN_i^j$ in Figure \ref{cones} are called the \textit{caps} of the nanocone. 

The following Corollary is another consequence of Theorem \ref{lowerbound}, and its proof, containing additional details, is provided in Section \ref{secnanocone}.

\begin{corollary} \label{lowernanocone}\.\\
	\begin{enumerate}
  \item[(a)] $s(CN_2^s(\ell)) \ge \frac{14+4\ell(\ell+4)}{3} - \frac{1}{3} - \frac{4(\ell+2)}{6}$
  \item[(b)] $s(CN_2^a(\ell)) \ge \frac{11+2\ell(2\ell+7)}{3} - \frac{1}{3	}-\frac{4\ell+7}{6}$
  \item[(c)] $s(CN_3^s(\ell)) \ge \frac{10+3\ell(\ell+4)}{3} - \frac{1}{2	}-\frac{3(\ell+2)}{6}$
  \item[(d)] $s(CN_3^a(\ell)) \ge \frac{16+\ell(3\ell+16)}{3} - \frac{1}{2	}-\frac{3\ell+8}{6}$
  \item[(e)] $s(CN_4^s(\ell)) \ge \frac{12+2\ell(\ell+6)}{3} - \frac{2}{3	}-\frac{2(\ell+3)}{6}$
  \item[(f)] $s(CN_4^a(\ell)) \ge \frac{15+2\ell(\ell+7)}{3} - \frac{2}{3	}-\frac{2\ell+7}{6}$
  \item[(g)] $s(CN_5^a(\ell)) \ge \frac{16+\ell(\ell+12)}{3} - \frac{5}{6	}-\frac{\ell+6}{6}$
	\end{enumerate}
\end{corollary}

The upper bound on the saturation number of nanocones relies on splitting the nanocone in subgraphs, where the number of subgraphs depends on the number of break edges. The following Lemma was proven in \cite{graves1}.

\begin{lemma} \label{break} \cite{graves1}
	In a nanocone, the number of pentagons $p$ and the number of break edges $s$ are related by 
	$$
	s+p=6.
	$$
\end{lemma}

\begin{comment}
A \textit{benzenoid parallelograms }, $P_{p,q}$, consists of a configuration of $p \times q$ congruent regular hexagons arranged in $p$ rows, each row containing $q$ hexagons, and each row shifted for half a hexagon to the right from the row immediately below.

\begin{proposition}\cite{doslic2} \label{benzparallel}
$$
s(P_{p,q}) \le \left \lceil \frac{2p+1}{3} \right \rceil q+p
$$
\end{proposition} 
\end{comment}

\begin{figure}[h]
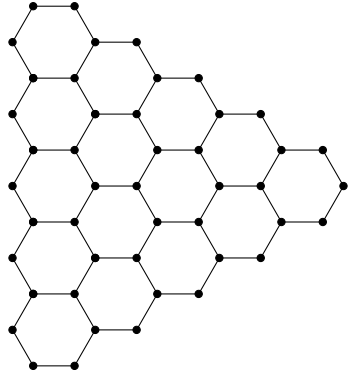

\begin{center}
	  \include{triangle}
\end{center}
  \caption{The benzenoid triangle, $T_5$.} \label{benztriangleT_5}
\end{figure}

\noindent The subgraph used is a \textit{benzenoid triangle}, $T_p$, which is a patch that can be constructed by arranging ${p+1 \choose 2}$ hexagonal faces in the shape of an equilateral triangle, so that each side of the triangle has $p$ hexagons. For an example of a benzenoid triangle, see Figure \ref{benztriangleT_5}. Note that the saturation number of similar graphs, such as benzenoid parallelograms, was studied in \cite{doslic2}. The upper bound on the saturation number of benzenoid triangles is presented in Lemma \ref{benztriangle} which is used to deduce the upper bounds on the saturation number of nanocones in Theorem \ref{uppernanocone}. The proofs of these results are presented in Section \ref{secnanocone}.

\begin{lemma} \label{benztriangle}
$$
  s(T_p) \le \left \lfloor \frac{(p+1)(p+3)}{3} \right \rfloor
$$
\end{lemma}

\noindent The upper bound presented in Lemma \ref{benztriangle} is believed to be the exact value of $s(T_p)$ and is the sequence A032765 in the OEIS \cite{oeis}.

\begin{theorem} \label{uppernanocone} \.\\
	\begin{enumerate}
	\item[(a)] $s(CNC_k(\ell)) \le k \left \lfloor \frac{\ell(\ell+2)}{3} \right \rfloor + k(\ell+1)$
  \item[(b)] $s(CN_2^s(\ell)) \le 4 \left \lfloor \frac{(\ell+1)(\ell+3)}{3} \right \rfloor + 4(\ell + 1) + 1$
  \item[(c)] $s(CN_2^a(\ell)) \le 3 \left \lfloor \frac{(\ell+1)(\ell+3)}{3} \right \rfloor + \left \lfloor \frac{\ell (\ell+2)}{3} \right \rfloor + 4(\ell + 1) + 1$
  \item[(d)] $s(CN_3^s(\ell)) \le 3 \left \lfloor \frac{(\ell+1)(\ell+3)}{3} \right \rfloor + 3(\ell + 1) + 1$
  \item[(e)] $s(CN_3^a(\ell)) \le 2 \left \lfloor \frac{(\ell+2)(\ell+4)}{3} \right \rfloor + \left \lfloor \frac{(\ell+1)(\ell+3)}{3} \right \rfloor + 3(\ell + 1) + 2$
  \item[(f)] $s(CN_4^s(\ell)) \le 2 \left \lfloor \frac{(\ell+2)(\ell+4)}{3} \right \rfloor + 2(\ell + 1) + 3$
  \item[(g)] $s(CN_4^a(\ell)) \le \left \lfloor \frac{(\ell+3)(\ell+5)}{3} \right \rfloor + \left \lfloor \frac{(\ell+2)(\ell+4)}{3} \right \rfloor + 2(\ell + 1) + 4$
  \item[(h)] $s(CN_5^a(\ell)) \le \left \lfloor \frac{(\ell+5)(\ell+7)}{3} \right \rfloor + (\ell + 1) + 6$
	\end{enumerate}
\end{theorem}

Combining Corollaries \ref{lowersd} and \ref{lowernanocone} along with Theorem \ref{uppernanocone} shows that if $G$ is any nanocone graph with $n$ vertices, then $s(G)\sim n/3$. Hence, in a smallest maximal matching, as $n$ gets large there are roughly 2 matched edges per hexagon. These findings are consistent with the work done on the saturation number of fullerenes \cite{andova} and benzenoid graphs \cite{doslic2}.

\begin{comment}
	\begin{corollary} \label{oldupsatcone}
$$
s(CNC_k[\ell]) \le \begin{cases} \frac{k}{2}\left ( \left \lceil \frac{2\ell+1}{3} \right \rceil \cdot \ell+2\ell+1 \right ) & \text{if } k \text{ is even}\\
\left \lfloor \frac{k}{2} \right \rfloor \left ( \left \lceil \frac{2\ell+1}{3} \right \rceil \cdot \ell + \ell \right ) + \frac{1}{2} \left( \left \lceil \frac{2\ell-1}{3} \right \rceil \cdot \ell +(\ell-1) \right) + \left \lceil \frac{k}{2} \right \rceil(\ell+1) & \text{if } k \text{ is odd} \end{cases}
$$
	\end{corollary}
\end{comment}

\subsection{Nanotubes}

\begin{figure}[h] 
	\begin{center}
		\input{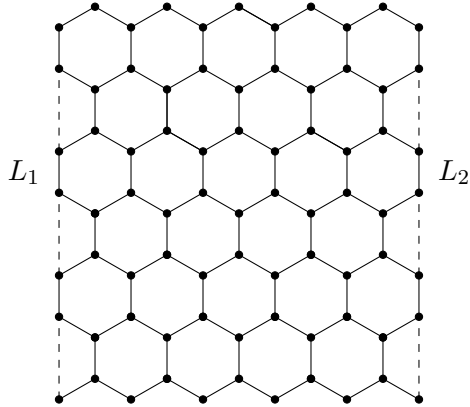}	
	\end{center}
	\caption{The zig-zag tubulene, $ZT(6,5)$, is obtained from the figure above by gluing the lines $L_1$ and $L_2$ together.} \label{zigzag}
\end{figure}

\begin{figure}[h]
	\begin{center}
		\input{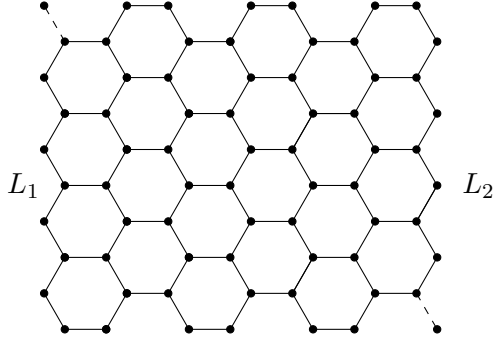}
	\end{center}
	\caption{The armchair tubulene, $AT(4,6)$, is obtained from the figure above by gluing the curves $L_1$ and $L_2$ together.} \label{armchair}
\end{figure}

\textit{Open ended nanotubes}, also called \textit{tubulenes}, can be obtained in the following way. Starting with a hexagonal tessellation of a cylinder, take the finite graph induced by all hexagons that lie between two vertex disjoint cycles, where each cycle encircles the axis of the cylinder. This paper considers two types of tubulenes having particularly nice structure, namely zig-zag and arm chair tubulenes, shown in Figures \ref{zigzag} and \ref{armchair}. 

Zig-zag tubulenes, $ZT(\ell,m)$, have $\ell$ horizontal layers of hexagons, each containing $m$ hexagons. Bounds on the saturation number of such zig-zag were first established in \cite{tratnik1}, shown in Corollary \ref{tratnikzigzag}.

\begin{corollary} \cite{tratnik1} \label{tratnikzigzag}
$$
\frac{m(\ell+1)}{2} \le s(ZT(\ell,m)) \le \begin{cases} \frac{m(2\ell+3)}{3} & \text{if } 3|\ell \\ \frac{m(2\ell+1)}{3} & \text{if } 3|(\ell-1) \\ \frac{m(2\ell+2)}{3} & \text{if } 3|(\ell-2) \end{cases}
$$
\end{corollary}

\noindent Corollary \ref{lowerzigzag} follows as an application of Theorem \ref{lowerbound} and improves the lower bound for the saturation number of zig-zag tubulenes. The proof is contained in Section \ref{sectubes}.

\begin{corollary} \label{lowerzigzag}
$$
s(ZT(\ell,m)) \ge \frac{m(2\ell+1)}{3}
$$
\end{corollary}

Combining the new lower bound from Corollary \ref{lowerzigzag} and the upper bounds presented in Corollary \ref{tratnikzigzag}, it follows that $s(ZT(\ell,m)) = \frac{m(2\ell+1)}{3}$ whenever $3|(\ell-1)$.

Armchair tubulenes, $AT(m,\ell)$, have $\ell$ vertical layers of hexagons, each containing $m$ hexagons. The saturation number of armchair tubulenes was also studied in \cite{tratnik1}, as seen in Corollary \ref{tratnikarmchair}.

\begin{corollary} \cite{tratnik1} \label{tratnikarmchair}
$$
\frac{\ell(m+1)}{2} \le s(AT(m,\ell)) \le \begin{cases} \frac{2\ell(m+1)}{3} & \text{if } 3|\ell \\ \frac{(2\ell+1)(m+1)}{3} & \text{if } 3|(\ell-1) \\ \frac{2(\ell+2)(m+1)}{3} & \text{if } 3|(\ell-2) \end{cases}
$$
\end{corollary}

\noindent Another application of Theorem \ref{lowerbound}, Corollary \ref{lowerarmchair} improves the lower bound for the saturation number of armchair tubulenes and its proof is found in Section \ref{sectubes}.

\begin{corollary} \label{lowerarmchair}
$$
s(AT(m,\ell)) \ge \frac{\ell(2m+1)}{3}
$$
\end{corollary}

From the work above, it follows that the saturation number of zigzag and armchair tubulenes with $n$ vertices is essentially $n/3$. This is now consistent with the findings for nanocones, fullerenes, and benzenoid graphs.

\section{Proof of the main tool}\label{seclowerbound}

\begin{proof}[Proof (of Theorem \ref{lowerbound})]
Let $M$ be a maximal matching in $G$. Let the edges in $M$ and the vertices saturated by $M$ be called black, and let the remaining edges and vertices be called white. Let $B$ and $W$ be the set of all black and white vertices, respectively.

The proof using the discharging method, setting the initial charges as follows:

\begin{itemize}
\item Let the initial charge of each black vertex be $3$;
\item Let the initial charge of each white vertex be $-6$;
\item Let the initial charge of each $k$-gonal, internal defect be equal to $\begin{cases}$k-2$ & \text{ if } k \text{ is odd} ,\\ k & \text{ if } k \text{ is even}\end{cases}$; and
\item Let the initial charge of each external defect, $E$, be $3n_2(E)$.
\end{itemize}

It remains to show that the total sum of the charge in the graph $3|B|-6|W|+\sum _{k \text{ odd}}(k-2)o_k+ \sum _{k \text{ even}} ke_k+3 n_2$ is non-negative. From this it follows that
$$
3|B|\ge 2|B|+2|W|-\frac{1}{3}\left( \sum _{k \text{ odd}}(k-2)o_k+ \sum _{k \text{ even}} ke_k \right)-n_2
$$
implying that
\begin{align*}
|M|=\frac{|B|}{2} &\ge \frac{|B|+|W|}{3}-\frac{1}{18} \left( \sum _{k \text{ odd}}(k-2)o_k+ \sum _{k \text{ even}} ke_k \right)-\frac{n_2}{6} \\
& = \frac{n}{3}-\frac{1}{18} \left( \sum _{k \text{ odd}}(k-2)o_k+ \sum _{k \text{ even}} ke_k \right)-\frac{n_2}{6}.
\end{align*}

\noindent The initial charge is distributed as follows:

\begin{enumerate}
\item[(R1)] Each external defect sends +3 charge to each incident, degree 2 white vertex.
\end{enumerate}

First note that all vertices of degree 2 in $G$ are incident to an external defect. There are a total of $n_2(E)$ vertices of degree 2 incident to an external defect $E$, not all of them white vertices, so after applying (R1) all white vertices of degree 2 in $G$ now have charge -3. 

\begin{enumerate}
\item[(R2)] Each white vertex distributes its negative charge evenly among the adjacent black vertices.
\end{enumerate}

Since $M$ is a maximal matching, $W$ is an independent set in $G$, so no 2 white vertices are adjacent. The white vertices of degree 2 are adjacent to exactly 2 black vertices and sends -1.5 charge to each adjacent black vertex. All other white vertices are adjacent to 3 black vertices and sends -2 charge to each adjacent black vertex. After applying (R2), all white vertices have charge 0.

Let $v$ be a black vertex. Since $v$ is saturated by $M$, $v$ is adjacent to at least one black vertex and hence, $v$ is adjacent to at most 2 white vertices. After receiving charge 0, -1.5, -2, -3, -3.5, or -4 from (R1) according to the number and type of white neighbors, $v$ now has charge 3, 1.5, 1, 0, -0.5, or -1. 

Next, let $e_v$ be the black edge incident with $v$, and let $f_v$ be the face incident to $v$ but not $e_v$, if such a face exists. Note that it's possible $f_v$ is an external defect. If such a face does not exists, then both $v$ and $e_v$ must be incident to an external defect, in which case set $u_v$ to be the incident external defect. 

\begin{enumerate}
\item[(R3)] Each black vertex sends all of its remaining charge to $f_v$ or $u_v$.
\end{enumerate}

Note that all charge that was initially present at the vertices of $G$ is now at its faces. Due to the face that $G$ is pseudoconvex, it is straightforward to check that if $v$ is a black vertex that sent charge to an external defect according to (R3), then $v$ previously had charge 0, 1, 1.5, or 3. So external defects receive no negative charge after applying (R3), and hence, their total charge is non-negative.

Now the only case when a face receives negative charge from (R3) is when a black vertex $v$ with 2 white neighbors sends charge $-1$ or $-\frac{1}{2}$ (depending on the degrees of the white neighbors) to $f_v$. So if a face is incident with at most 1 white vertex, then its charge is certainly non-negative. It turns out that if a face is incident to at most two white vertices, then its charge is non-negative.

An internal $k$-gonal defect is incident to at most $k/2$ white vertices if $k$ is even and $(k-1)/2$ white vertices if $k$ is odd. Hence, the negative charge an internal $k$-gonal defect receives after applying (R3) is at most $k/2$ in either case. Therefore, the charge of each such defect is at least $\lfloor k/2 \rfloor$ after applying (R3).

All hexagonal faces have nonnegative charge except those incident to three white vertices. For the hexagons incident to at least two white vertices, the different cases for hexagons are split into figures depending on the number of incident vertices of degree 2. The cases for hexagons incident to 0, 1, 2, and 3 vertices of degree 2 can be seen in Figures \ref{internalhex}, \ref{externalhex1}, \ref{externalhex2}, and \ref{externalhex3}, respectively. 

\begin{figure}[h]
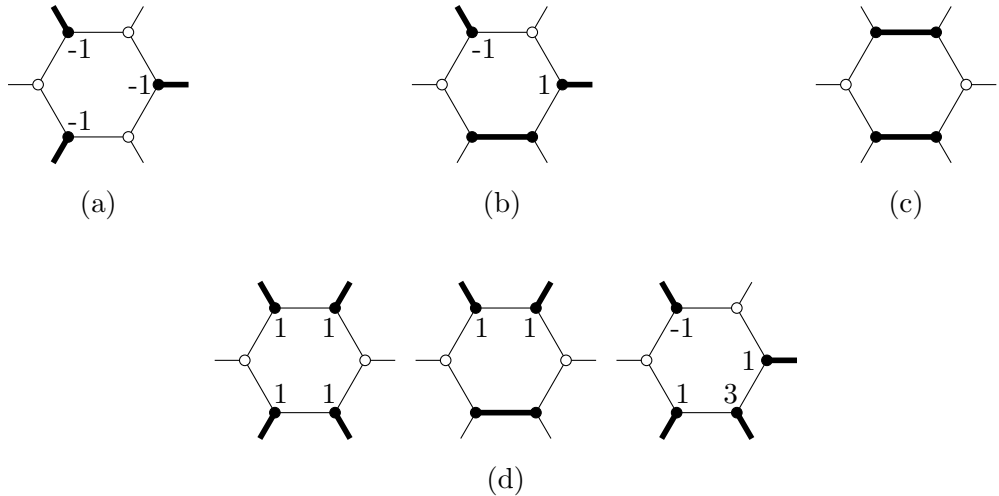
 
\begin{center}
\input{h-bad} \hspace{1in} \input{h-transition} \hspace{1in} \input{h-neutral} 

\vspace{0.25in}

\input{h-good1} \input{h-good2} \input{h-good3}
\end{center}
\caption{Hexagons adjacent to at least 2 white vertices and 0 vertices of degree 2.} \label{internalhex}
\end{figure}

\begin{figure}[h]
\begin{center}
\input{he-bad} \hspace{1in} \input{he-transition1} \hspace{-0.25in} \input{he-transition-sep} \hspace{-0.25in} \input{he-transition2} 

\vspace{0.25in}

\input{he-neutral1} \hspace{-0.25in} \input{he-neutral-sep} \hspace{-0.25in} \input{he-neutral2} \hspace{1in} \input{he-nuetral3}

\vspace{0.25in}

\input{he-good1} \input{he-good2} \input{he-good3} \input{he-good4} \input{he-good5}
\end{center}
\caption{Hexagons adjacent to at least 2 white vertices and 1 vertex of degree 2.} \label{externalhex1}
\end{figure}

\begin{figure}[h]
\begin{center}
\input{he2-transition} \hspace{1in} \input{he2-neutral1} \hspace{-0.25in} \input{he2-neutral-sep} \hspace{-0.25in} \input{he2-neutral2}

\vspace{0.25in}

 \input{he2-neutral3} \hspace{1in} \input{he2-good1} \hspace{-0.25in} \input{he2-good-sep} \hspace{-0.25in} \input{he2-good2}

\end{center}
\caption{Hexagons adjacent to at least 2 white vertices and 2 vertices of degree 2.} \label{externalhex2}
\end{figure}

\begin{figure}[h]
\begin{center}
\input{he3-neutral1} \hspace{-0.25in} \input{he3-neutral-sep} \hspace{-0.25in} \input{he3-neutral2} \hspace{1in} \input{he3-neutral3} \hspace{1in} \input{he3-good1}

\end{center}
\caption{Hexagons adjacent to at least 2 white vertices and 3 vertices of degree 2.} \label{externalhex3}
\end{figure}

If a hexagon has negative charge as in Figure \ref{internalhex} (a) or Figure \ref{externalhex1} (a), then these hexagons are called \textit{bad}. If a hexagon has zero charge as in Figure \ref{internalhex} (b), Figure \ref{externalhex1} (b), or Figure \ref{externalhex2} (a), then these hexagons are called \textit{transitional}. If a hexagon has zero charge as in Figure \ref{internalhex} (c), Figure \ref{externalhex1} (c), Figure \ref{externalhex2} (b), or Figure \ref{externalhex3} (a), then these hexagons are called \textit{neutral}. Those hexagons with charge $1/2$ as in Figure \ref{externalhex1} (d), Figure \ref{externalhex2} (b), or Figure \ref{externalhex3} (b) are called \textit{almost neutral}. All other hexagons have a positive charge, and the value of the positive charge is at least the number of incident white vertices. These hexagons with positive charge are called \textit{good}.

Let $f$ be a transitional hexagon. Then $f$ is incident to one black edge, two white vertices, and two black vertices that are incident to black edges not incident to $f$. Let the white vertex adjacent to the black edge incident to $f$ be called \textit{outgoing}. Let the other white vertex, between the black edges that are not incident to $f$, be called \textit{incoming}.

The last steps of the discharging are given by the following rules:

\begin{enumerate}
\item[(R4)] Each good face sends charge 1 to each incident, degree 3 white vertex.

\item[(R5)] Each bad hexagonal face sends charge -1 to each incident, degree 3 white vertex.

\item[(R6)] Each transitional hexagonal face sends charge -1 to the incoming degree 3 white vertex, and it sends charge 1 to the outoing degree 3 white vertex.

\end{enumerate}

After applying (R4)-(R6), there is no negative charge left at the faces of $G$, and the only possible negative charge resides at white vertices of degree 3 in $G$. 

Let $v$ be a vertex that was sent charge $-1$ by either (R5) or (R6), let $h$ be the hexagon that sent the negative charge to $v$, and let $u_i$, $i=1,2$, be the black vertices adjacent to $v$ and incident to $h$. Since $h$ is either a bad hexagon or transitional hexagon, then the black edges incident to the $u_i$ are not incident to $h$. 

Now let $x$ be the black vertex adjacent to $v$ but not incident to $h$. Since $G$ is pseudoconvex, there exists two faces $f_i$, $=1,2$, incident with $v$ different from $h$, where $f_i$ is incident to $u_i$, $i=1,2$. Without loss of generality, assume the black edge incident with $x$ is incident with $f_1$. Since the black edges incident to both $u_1$ and $x$ are incident to $f_1$, it follows that $f_1$ is either good or neutral, so it does not send negative charge to $v$.

Now consider $f_2$. Since the black edge incident to $u_2$ is incident to $f_2$, then $f_2$ cannot be a bad hexagon, nor a neutral hexagon due to the black edge incident to $x$. Furthermore, $f_2$ cannot be almost-neutral since $v$ is a degree 3 white vertex. If $f_2$ is not incident to any other white vertex other than $v$, then $f_2$ is a good hexagon. If $f_2$ is incident to another white vertex at distance 3 from $v$, then it is a good hexagon as well. If $f_2$ is incident to another white vertex at distance 2 from $v$, then $f_2$ is a transitional hexagon. In the case that $f_2$ is transitional, then $v$ is the outgoing white vertex for $f_2$. Hence in all considered cases, $f_2$ has sent positive 1 charge to $v$ by (R4) or (R5).

\begin{figure}[h]
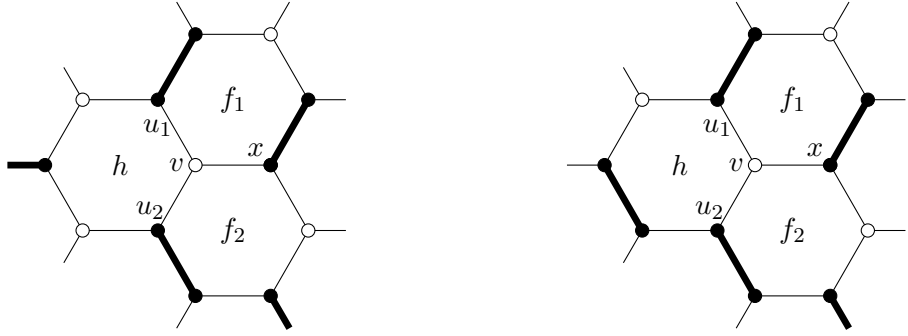

\begin{center}
\input{f2bad} \hspace{1in} \input{f2transition}
\end{center}
\caption{On the left shows the case when $h$ is a bad hexagon and $f_2$ is transitional. On the right shows when both $h$ and $f_2$ are transitional.}\label{f2}
\end{figure}

Repeating the above argument, it's possible $G$ has a chain of adjacent, transitional hexagons, which in turn, would move charge between adjacent hexagons. If such a chain starts with a bad hexagon, then the chain cannot close on itself forming a cycle of hexagons. Such a cycle would have to close at the bad hexagon, implying a white vertex receives negative charge from both a bad hexagon and transitional hexagon and this cannot happen according to the above argument. A chain beginning with a transitional hexagon could close to form a cycle of transitional hexagons, in which case, since transitional hexagons have zero charge, the discharging simply moved zero charge around in a cycle.

Thus after these last steps of discharging, there is no negative charge in the graph. So the total sum of charge is non-negative, which finishes the proof.

\end{proof}

\section{Proofs for nanocones} \label{secnanocone}

\begin{proof}[Proof (of Lemma \ref{benztriangle})]
	First a construction of a maximal matching $M$ is given and then below it is shown this matching yields the desired bound. To construct $M$, $T_p$ is drawn in the plane so that the hexagons appear in columns and the number of hexagons in columns decreases moving to the right, as in Figure \ref{benztrianglematch}. Moving left to right, the following pattern of matched edges is iterated every 3 columns of hexagons: the first column of $k$ hexagons requires $k+1$ matched edges, the second column of $k-1$ hexagons requires $k$ matched edges, and the third column is skipped, since edges from the second column partially matches the third column. See the bold edges in Figure \ref{benztrianglematch} for examples of these matchings. This pattern is continued so long as there are at least 3 columns of hexagons remaining, at which point the pattern breaks. This process yields a maximal matching of size
	
	$$
	(p+1)+p+(p-2)+(p-3)+ \cdots + k_2 + k_1
	$$
	
	where the end values $k_i$, $i=1,2$, fall into 3 cases depending on the value of $p+1$ modulo 3, which can be seen in the matchings of $T_3$, $T_4$, and $T_5$ in Figure \ref{benztrianglematch}. 
	
	\begin{figure}[h]
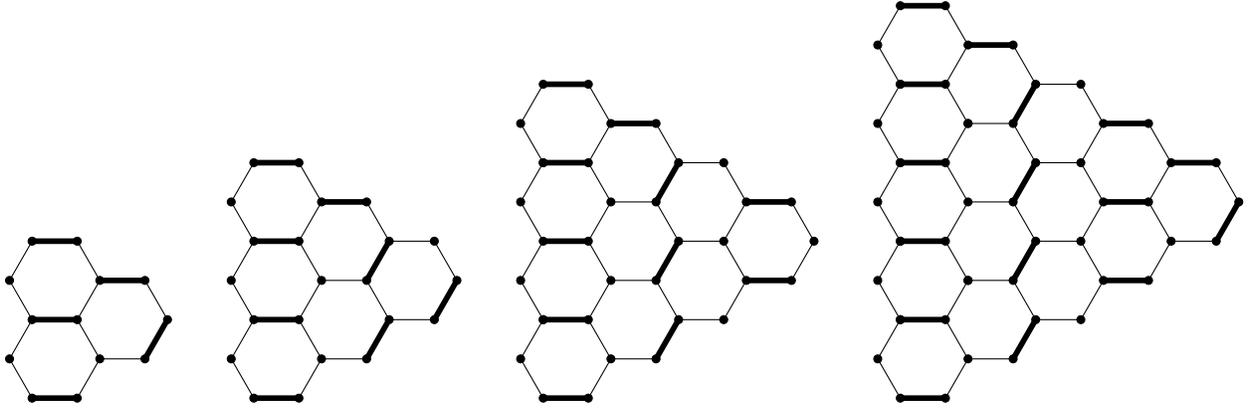

	\input{triangle2} \hfill \input{triangle3} \hfill \input{triangle4} \hfill \input{triangle5}
	
	\caption{Maximal matchings of $T_2$, $T_3$, $T_4$, and $T_5$ as described in Lemma \ref{benztriangle}.} \label{benztrianglematch}
	\end{figure}
	
	The remaining argument is broken into these 3 cases:
	
	\begin{enumerate}
  \item[(Case 1)] $p+1 \equiv 0 \pmod 3$
	\end{enumerate}
	
	In this case, the construction gives a matching of size
	$$
	((p+1)+p) + ((p-2)+(p-3)) + \cdots + (6+5) + (3+2)  
	$$
	and then iteratively rearranging terms so that the next largest term is now paired with the next smallest term yields
	$$
	((p+1)+2) + (p+3) + ((p-2)+5) + ((p-3)+6) + \cdots ((p-k)+(k+3))  
	$$
	for some value $k$. This new sum consists of $(p+1)/3)$ pairs each summing to $(p+3)$, so the matching has size exactly $(p+1)(p+3)/3$.
	
	\begin{enumerate}
  \item[(Case 2)] $p+1 \equiv 1 \pmod 3$
	\end{enumerate}
	
	In this case, it also follows that $(p+3) \equiv 0 \pmod 3$. The construction gives a matching of size
	$$
	((p+1)+p) + ((p-2)+(p-3)) + \cdots + (4+3) + 1  
	$$
	which can be rearranged to
	$$
	((p+1)+0) + (p+1) + ((p-2)+3) + ((p-3)+4) + \cdots ((p-k)+(k+1))  
	$$
	for some value $k$. This new sum consists of $(n+3)/3)$ pairs each summing to $(p+1)$, so the matching has size exactly $(p+1)(p+3)/3$.

	\begin{enumerate}
  \item[(Case 3)] $p+1 \equiv 2 \pmod 3$
	\end{enumerate}
	
	This case has $(p+1) = 3q + 2$ for some integer $q$. The construction gives a matching of size
	$$
	((p+1)+p) + ((p-2)+(p-3)) + \cdots + (5+4) + 2.  
	$$
	First, the smallest and largest terms are paired together, $((p+1)+2)$, and  the next largest term, $p$, is reserved. The remaining terms are iteratively rearranged so that the next largest term is now paired with the next smallest term to obtain
	$$
	((p-2)+4) + ((p-3)+5) + \cdots ((p-k)+(k+2))  
	$$
	for some value $k$. This last sum results in $(q-1)$ pairs each summing to $(p+2)$. Hence the matching has size $(p+3)+p+(q-1)(p+2)$. Now 
\begin{align*}
	(p+3)+p+(q-1)(p+2) &= (p+3)+p+(q-1)(p+3)-(q-1) \\
	&= p+q(p+3)-(q-1) \\
	&= p-\frac{p-4}{3} + q(p+3) \\
	&< \frac{2}{3}(p+3) + q(p+3) \\
	&=\frac{(p+1)(p+3)}{3}
\end{align*}
which proves the desired bound.

\end{proof}

\begin{proof}[Proof (of Corollary \ref{lowernanocone})]
	The claimed lower bounds are a straight forward application of Theorem \ref{lowerbound} depending on the total number of vertices, the number of pentagons, and the number of vertices of degree 2. The counts of these values are provided below, where both counts of vertices follow by induction on $\ell$.
		
		(a) $CN_2^s(\ell)$ has $14+4\ell(\ell+4)$ total vertices, 2 pentagons, and $4(\ell +2)$ vertices of degree 2.

		(b) $CN_2^a(\ell)$ has $11+2\ell(2\ell+7)$ total vertices, 2 pentagons, and $4\ell +7$ vertices of degree 2.
		
		(c) $CN_3^s(\ell)$ has $10+3\ell(\ell+4)$ total vertices, 3 pentagons, and $3(\ell +2)$ vertices of degree 2.
		
		(d) $CN_3^a(\ell)$ has $16+\ell(3\ell+16)$ total vertices, 3 pentagons, and $3\ell +8$ vertices of degree 2.
		
		(e) $CN_4^s(\ell)$ has $12+\ell(2\ell+12)$ total vertices, 4 pentagons, and $2(\ell+3)$ vertices of degree 2.
		
		(f) $CN_4^a(\ell)$ has $15+\ell(2\ell+14)$ total vertices, 4 pentagons, and $2\ell+7$ vertices of degree 2.
		
		(g) $CN_5^a(\ell)$ has $16+\ell(\ell+12)$ total vertices, 5 pentagons, and $\ell+6$ vertices of degree 2.
\end{proof}

\begin{figure}
\centering
	\begin{subfigure}[h]{.3\textwidth}
		\includegraphics[width=\textwidth]{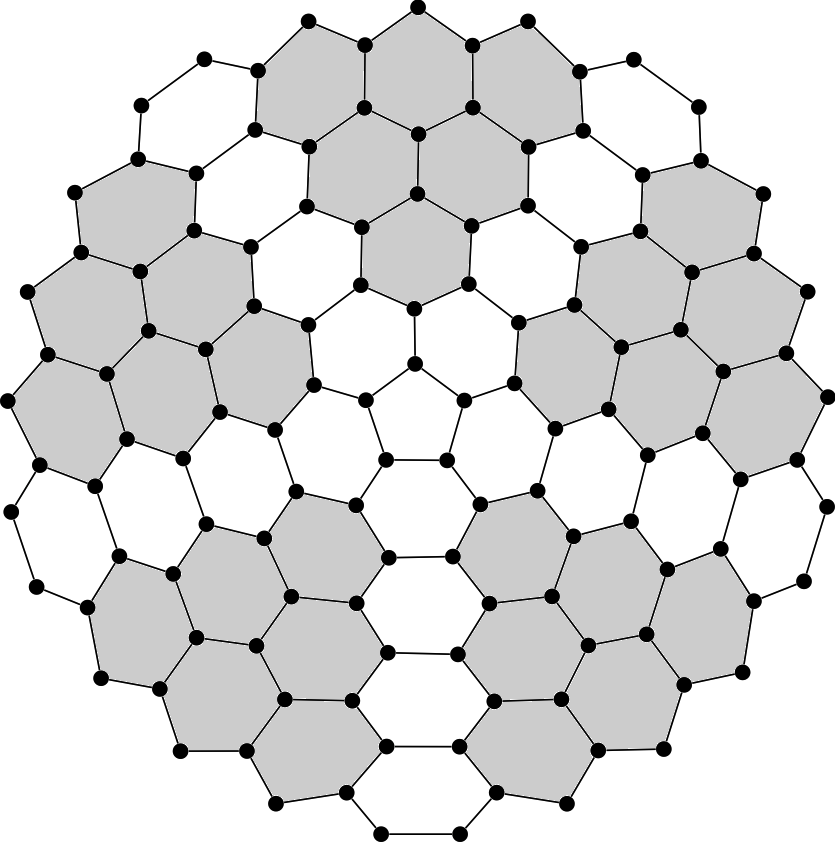} 
		\caption*{$CNC_5(4)$}
	\end{subfigure}
	\begin{subfigure}[h]{.3\textwidth}
		\includegraphics[width=\textwidth]{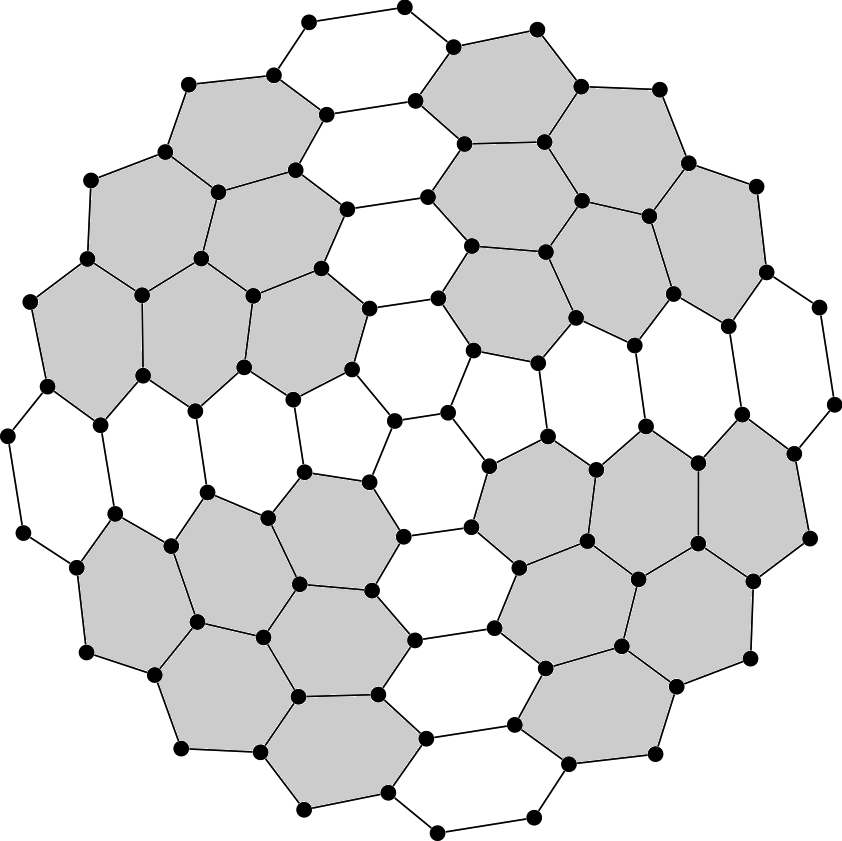} 
		\caption*{$CN_2^s(3)$}
	\end{subfigure}
	\begin{subfigure}[h]{.3\textwidth}
		\includegraphics[width=\textwidth]{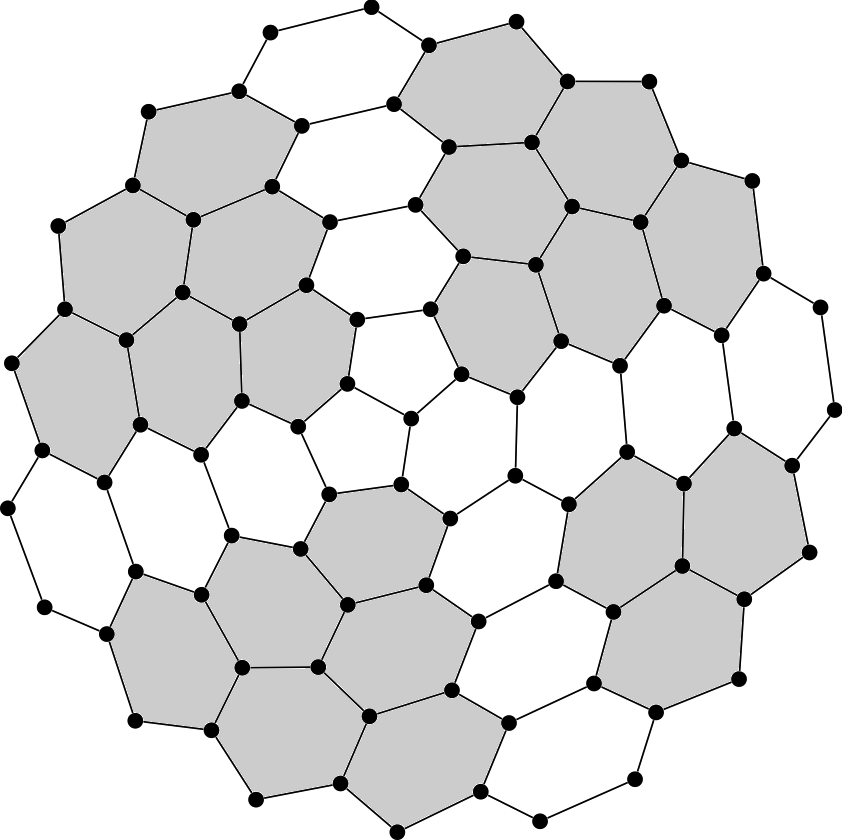} 
		\caption*{$CN_2^a(3)$}
	\end{subfigure}
	
	\begin{subfigure}[h]{.3\textwidth}
		\includegraphics[width=\textwidth]{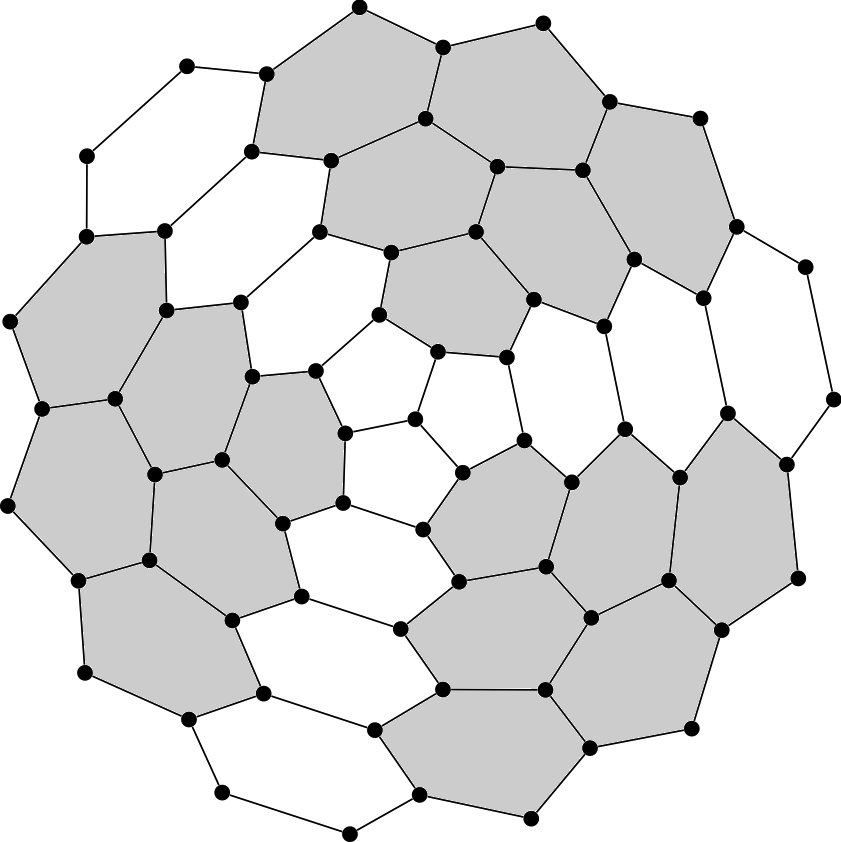} 
		\caption*{$CN_3^s(3)$}
	\end{subfigure}
	\begin{subfigure}[h]{.3\textwidth}
		\includegraphics[width=\textwidth]{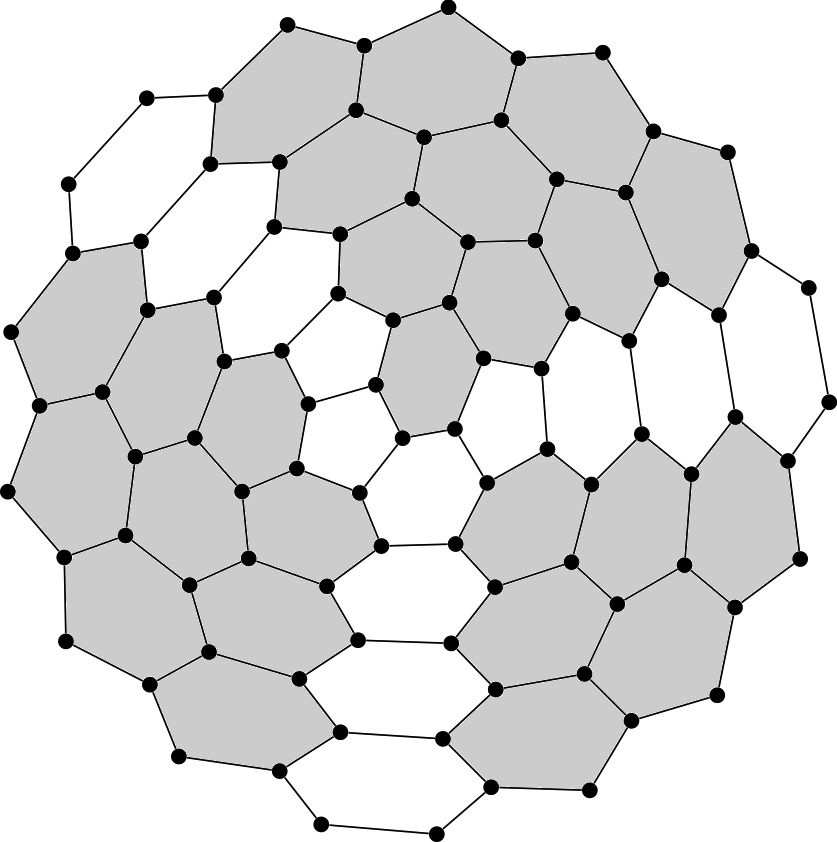} 
		\caption*{$CN_3^a(3)$}
	\end{subfigure}
	\begin{subfigure}[h]{.3\textwidth}
		\includegraphics[width=\textwidth]{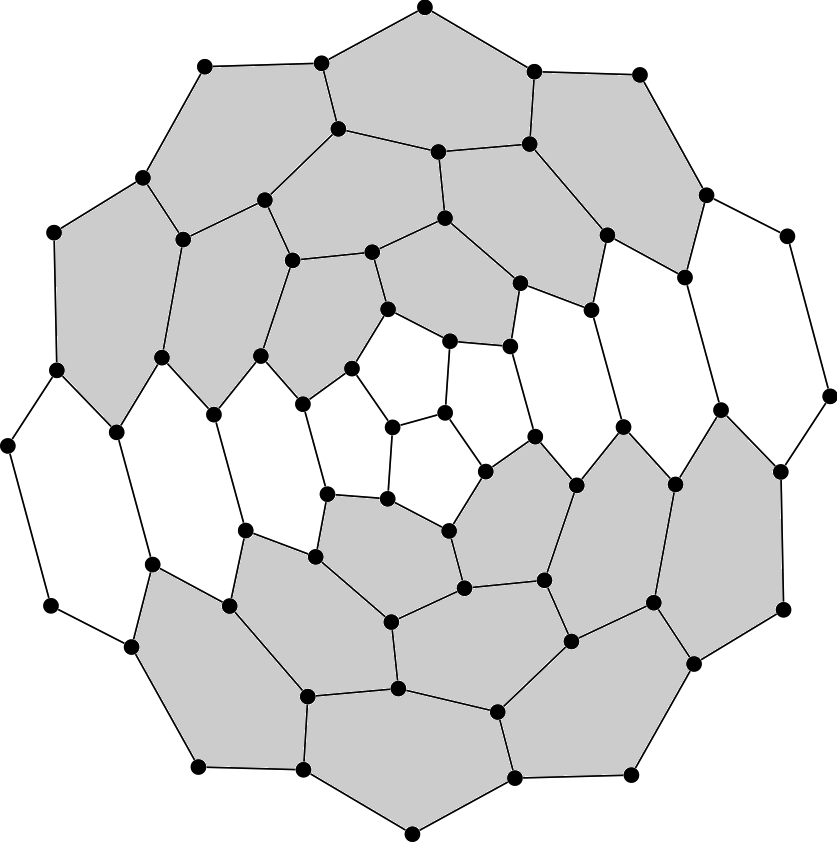}
		\caption*{$CN_4^s(3)$}
	\end{subfigure}
	
	\begin{subfigure}[h]{.3\textwidth}
		\includegraphics[width=\textwidth]{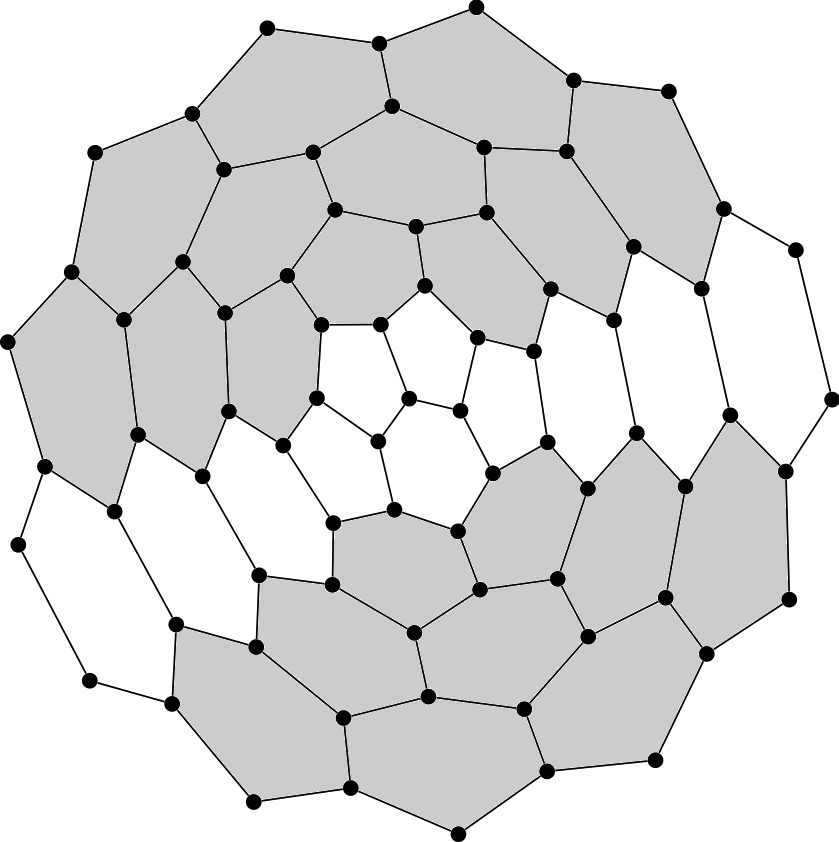}
		\caption*{$CN_4^a(3)$}
	\end{subfigure}
	\begin{subfigure}[h]{.3\textwidth}
		\includegraphics[width=\textwidth]{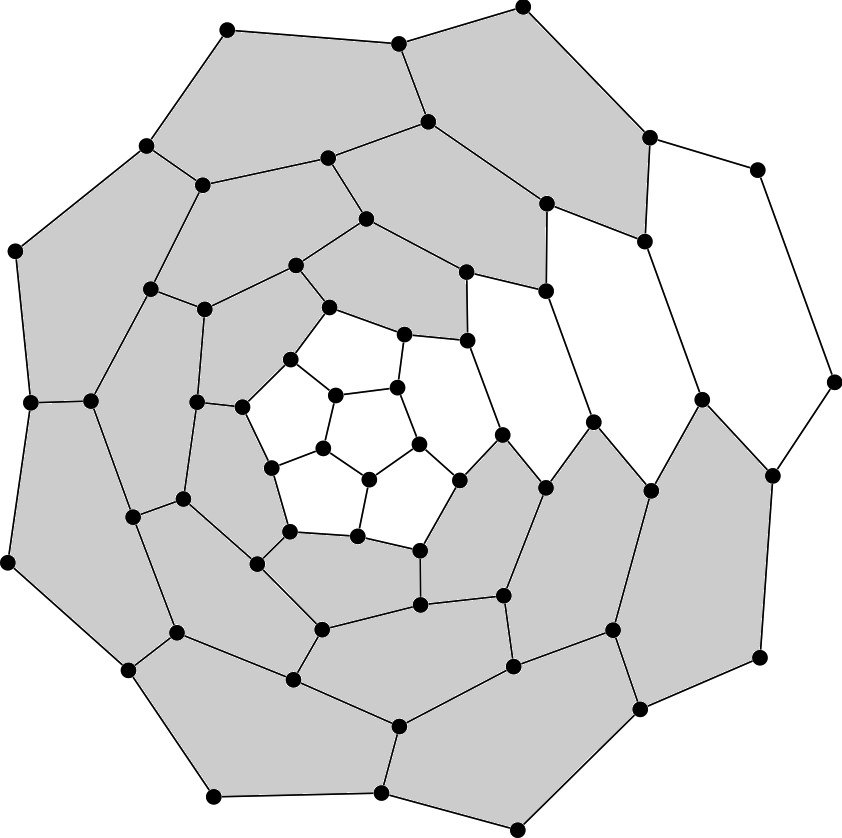}
		\caption*{$CN_5^a(3)$}
	\end{subfigure}
 
	\caption{Nanocones split into benzenoid triangles, or subgraphs thereof, which are represented by the shaded regions.} \label{nanoconetriangle}
\end{figure}

\begin{proof}[Proof (of Theorem \ref{uppernanocone})]
	Observe that a nanocone with $s$ break edges can be split into $s$ benzenoid triangles, or subgraphs of benzenoid triangles. Each such benzenoid triangle or subgraph resides between successive hexagons containing the break edges each layer of hexagons, see Figure \ref{nanoconetriangle}. The sizes of the triangles depends on the lengths of the sides of the nanocone. Lemma \ref{benztriangle} can be used to find a maximal matching of the benzenoid triangles of the indicated size. The union of these matchings augmented by a matching of size at most $s(\ell+1)$ along the break edges from each layer, and potentially an additional matching of the cap, gives an upper bound on the size of a maximal matching of the nanocone. Additional details for each case are provided below.
		
		(a) $CNC_k(\ell)$  has $k$ break edges and therefore can be split into $k$ benzenoid triangles $T_{\ell-1}$, each triangle with a matching of size $\left \lfloor \frac{\ell(\ell+2)}{3} \right \rfloor$. The union of these matching augmented by a matching of size of size $k(\ell+1)$ along the break edges from each layer provides an upper bound for a maximal matching of $CNC_k(\ell)$.
		
		(b) By Lemma \ref{break}, $CN_2^s(\ell)$ has 4 break edges and can be split into 4 benzenoid triangles, $T_\ell$, each with a matching of size $\left \lfloor \frac{(\ell+1)(\ell+3)}{3} \right \rfloor$. Their union augmented by a matching along the break edges of size at most $4(\ell +1)$ plus an addition edge needed for the remaining edges on the cap yields the desired upper bound.
		
		(c) By Lemma \ref{break}, $CN_2^a(\ell)$ has 4 break edges and can be split into 3 $T_\ell$'s each with a matching of size $\left \lfloor \frac{(\ell+1)(\ell+3)}{3} \right \rfloor$ and an additional $T_{\ell-1}$ with a matching of size $\left \lfloor \frac{\ell(\ell+2)}{3} \right \rfloor$. The break edges require at most $4(\ell +1)$ edges after which the cap requires 1 additional edge.
		
		(d) By Lemma \ref{break}, $CN_3^s(\ell)$ has 3 break edges and can be split into 3 $T_\ell$'s each with a matching of size $\left \lfloor \frac{(\ell+1)(\ell+3)}{3} \right \rfloor$. The break edges union the cap of the nanocone require at most an additional $3(\ell + 1)+1$ edges.
		
		(e) Again using Lemma \ref{break}, $CN_3^a(\ell)$ has 3 break edges and can be split into $T_{\ell + 1}$, $T_\ell$, and a subgraph of $T_{\ell+1}$, which in total require at most 
		$$
		2 \left \lfloor \frac{(\ell+2)(\ell+4)}{3} \right \rfloor + \left \lfloor \frac{(\ell+1)(\ell+3)}{3} \right \rfloor
		$$
		matched edges. The break edges need at most $3(\ell + 1)$ edges and the cap requires at most 2 edges, proving the desired bound. 
		
		(f) Lemma \ref{break} gives that $CN_4^s(\ell)$ has 2 break edges. So $CN_4^s(\ell)$ can be split into two pieces which turn out to be subgraphs of $T_{\ell+1}$, and each subgraph has a maximal matching of size at most $\left \lfloor \frac{(\ell+2)(\ell+4)}{3} \right \rfloor$. The union of these matchings augmented by a matching of the break edges of size at most $2(\ell+1)$ along with a matching of size 3 for the remaining edges of the cap provides the desired maximal matching.
		
		(g) Similar to the case in (f), $CN_4^a(\ell)$ can be split into subgraphs of $T_{\ell+1}$ and $T_{\ell+2}$ requiring at most 
		$$
		\left \lfloor \frac{(\ell+3)(\ell+5)}{3} \right \rfloor + \left \lfloor \frac{(\ell+2)(\ell+4)}{3} \right \rfloor
		$$ 
		matched edges. The break edges again require at most $2(\ell + 1)$ matched edges, after which the cap needs at most 4 edges.
		
		(h) By Lemma \ref{break}, $CN_5^a(\ell)$ has 1 break edge and contains a subgraph of $T_{\ell+4}$, which according to Lemma \ref{benztriangle} has a maximal matching of size at most $\left \lfloor \frac{(\ell+5)(\ell+8)}{3} \right \rfloor$. The break edges require at most $(\ell+1)$ matched edges and the cap needing an additional 6 matched edges.
		
		\end{proof}

	\section{Proofs for nanotubes} \label{sectubes}

\begin{proof}[Proof (of Corollary \ref{lowerzigzag})]
	It follows that $ZT(\ell, m)$ has $2m\ell + 2m$ total vertices and two external defects at the ends of the cylinder. The external defects each have $m$ vertices of degree 2, for a total of $2m$ degree 2 vertices. Now Theorem \ref{lowerbound} gives that
	\begin{align*}
		s(ZT(\ell, m)) &\ge \frac{2m\ell+2m}{3} - \frac{2m}{6} \\
		&= \frac{m(2\ell+1)}{3}.
	\end{align*}

\end{proof}

\begin{proof}[Proof (of Corollary \ref{lowerarmchair})]
	The armchair tubulene $AT(m, \ell)$ has $2m\ell + 2\ell$ total vertices and $2\ell$ vertices of degree 2. Theorem \ref{lowerbound} gives that
	\begin{align*}
		s(AT(m, \ell)) &\ge \frac{2m\ell+2\ell}{3} - \frac{2\ell}{6} \\
		&= \frac{\ell(2m+1)}{3}.
	\end{align*}

\end{proof}

\section{Acknowledgements}

This work was partially supported by a SPARC Graduate Research Grant from the Office of the Vice President for Research at the University of South Carolina. I am extremely grateful to Tomislav Do\v{s}li\'{c}, whose mentoring and hosting during the Summer 2015 eventually led to the completion of this work. I would also like  acknowledge the software CaGe \cite{cage} which helped draw Figures \ref{cones} and \ref{nanoconetriangle} in this paper, and I would like to thank Nico Van Cleemput, who helped export graphs generated with CaGe for my initial experimentation. I would like to thank Michael Santana for a helpful discussion concerning Theorem \ref{lowerbound}.

\bibliographystyle{amsplain}
\bibliography{ref}

%%%%%%%%%%%%%%%%%%%%%%%%%%%%%%%%%%%%%%%%%%%%%%%%%%%%
\end{document}